\documentclass[12pt,leqno]{article}

\usepackage{lineno}
\modulolinenumbers[5]
\usepackage{caption}
\usepackage[a4paper]{geometry}
\usepackage{amsfonts}
\usepackage{arydshln}
\usepackage{xtab}
\usepackage{setspace}
\usepackage{fancyhdr}
\usepackage{graphicx}
\usepackage{subcaption}
\usepackage{epstopdf}
\usepackage{resizegather}
\usepackage{mathtools}
\usepackage{hyperref}
\usepackage[T1]{fontenc}
\usepackage[latin1,utf8]{inputenc}
\usepackage[main=english,italian]{babel}

\usepackage{sidecap}

\sidecaptionvpos{figure}{c}

\usepackage{type1ec}
\usepackage[final]{microtype}
\usepackage{amsmath,amssymb,amsthm,eucal,dsfont,bm,mathrsfs,stmaryrd}
\usepackage{lmodern}
\usepackage{textcomp}
\usepackage{pict2e,enumitem}
\usepackage{hyperref}
\usepackage[dvipsnames]{xcolor}
\usepackage[nodayofweek]{datetime} 
\usepackage{comment}

\definecolor{cornellred}{rgb}{0.7, 0.11, 0.11}

\hypersetup{
    pdfpagemode={UseOutlines},
    bookmarksopen,
    pdfstartview={FitH},
    colorlinks,
    linkcolor={black},
    citecolor={cornellred},
    urlcolor={ForestGreen}
            }

\definecolor{cornellred}{rgb}{0.7, 0.11, 0.11}

\usepackage{geometry}

\newdateformat{monthyear}{\monthname[\THEMONTH], \THEYEAR}

\def\R0{\mathcal{R}_0}

\theoremstyle{plain}
        \newtheorem{theorem}{Theorem}[section]

\theoremstyle{definition}

\renewcommand
        {\thefootnote}{\arabic{footnote}}
        
\newcommand{\symfootnote}[1]{%
\let\oldthefootnote=\thefootnote%
\stepcounter{mpfootnote}%
\addtocounter{footnote}{-1}%
\renewcommand{\thefootnote}{\fnsymbol{mpfootnote}}%
\footnote{#1}%
\let\thefootnote=\oldthefootnote%
}

\usepackage{tikz}
\tikzstyle{comp}=[circle,draw=black!50,thick,minimum size=1cm]
\usetikzlibrary{shapes.geometric}
\usepackage{floatrow}

\counterwithout{equation}{section}

\title{A survey on Lyapunov functions for epidemic compartmental models}
\author{N. Cangiotti\footnote{Politecnico di Milano, Department of Mathematics, via Bonardi 9, Campus Leonardo, 20133, Milan (Italy). E-mail: \texttt{nicolo.cangiotti@polimi.it}}, \ M. Capolli\footnote{Institute of Mathematics, Polish Academy of Sciences, Jana i Jedrzeja Sniadeckich 8, 00-656, Warsaw, (Poland). E-mail: \texttt{mcapolli@impan.pl}}, \ M. Sensi\footnote{MathNeuro Team, Inria at Universit\'e C\^ote d'Azur, 2004 Rte des Lucioles, 06410, Biot, (France). E-mail: \texttt{mattia.sensi@inria.fr}}, \ S. Sottile\footnote{Department of Mathematics, University of Trento, Via Sommarive 14, 38123, Povo, (Italy). E-mail: \texttt{sara.sottile@unitn.it}}}
\date{}

\begin{document}
\mathtoolsset{showonlyrefs=true}

\maketitle

\begin{abstract}
In this survey, we propose an overview on Lyapunov functions for a variety of compartmental models in epidemiology. We exhibit the most widely employed functions, together with a commentary on their use. Our aim is to provide a comprehensive starting point to readers who are attempting to prove global stability of systems of ODEs. The focus is on mathematical epidemiology, however some of the functions and strategies presented in this paper can be adapted to a wider variety of models, such as prey-predator or rumor spreading.
\medskip

\paragraph{Mathematics Subject Classification:} 34D20, 34D23, 37N25, 92D30.
\medskip

\paragraph{Keywords:} Epidemic models, Lyapunov functions, Compartmental models, Global stability, Ordinary Differential Equations, Disease Free and Endemic Equilibria.
\end{abstract}

\section{Introduction}\label{sec1}

Stemming from the pioneering work of Kermack and McKendrick \cite{kermack1927contribution}, the mathematical modelling of infectious diseases has developed, over the last century, in various directions. An abundance of approaches and mathematical techniques have been employed to capture the many facets and details which describe the spread of an infectious disease in a population.

In particular, compartmental models remain one of the most widely employed approaches. In these models, a population is partitioned into compartments, characterizing each individual with respect to its current state in the epidemic. One can then write a system of Ordinary Differential Equations (from here onwards, ODEs) to study the evolution in time of the disease.

These models usually take their names from the compartments they consider, the most renowned one being the Susceptible-Infected-Recovered (SIR) model. The SIR models can be extended to SIRS models by considering the acquired immunity to be temporary rather than permanent, allowing Recovered individuals to become Susceptible again. Various compartments can be added, depending on the characteristic of the specific disease under study: Asymptomatic, Exposed, Waning immunity and many others.
\medskip

A remarkably useful tool for the study of this kind of models are Lyapunov functions, which ensure global (or, in some cases, local) asymptotic convergence towards one of the equilibria of the system. 

Given a system of $n$ ODEs $X'=f(X)$ and an equilibrium point $X^*$, we call a  scalar function $V\in C^1(\mathbb{R}^n,\mathbb{R})$ a \emph{Lyapunov function} if the following hold:
\begin{enumerate}
    \item $V$ attains its minimum at $X=X^*$;
    \item $V'=\nabla V \cdot f<0$ for $X\neq X^*$.
\end{enumerate} The classical definition of Lyapunov function requires also the conditions
\begin{enumerate}
    \item[3.] $X^*=0$ and $V(X^*)=0$;  
\end{enumerate}however, these amount to a change of coordinates in $\mathbb{
R}^n$ and a vertical translation of $V$, so we will accept the more general definition. The existence of such a function guarantess the global stability of the equilibrium $X^*$, as orbits of the systems naturally evolve towards the minimum power level of $V$.

The Basic Reproduction Number $\R0$ is a well know threshold in epidemics models. Usually, $\R0<1$ suggests Global Asymptotic Stability (from here onwards, GAS) of the Disease Free Equilibrium (from here onwards, DFE), whereas $\R0>1$ suggests GAS of the Endemic Equilibrium (from here onwards, EE).
In more complex models, the aforementioned conditions on $\R0$ might not be sufficient to prove the GAS of either equilibria, especially in cases in which the EE is not unique.
Lyapunov functions often explicitly involve $\R0$ to guarantee the extinction of the disease or its endemicity over time.

Unfortunately, given a generic system of ODEs, there is no universal way of deriving a Lyapunov function, nor to rule out the existence of one. However, there exist a few Lyapunov functions which have proven quite effective in a variety of different models.
\medskip

In this survey, we collect some of the most relevant functions available in the literature, to provide the reader with a series of options to apply to the model of their interest, depending on its formulation. We include an extensive bibliography to complement the essential information of each model we present. This will provide the reader with a convenient starting point to investigate the availability of a known Lyapunov function to analytically prove the asymptotic behaviour of their system of ODEs. For the sake of brevity, we do not repeat the proofs to show that any of the functions we present are, indeed, Lyapunov function for the respective system of ODEs. These proofs can be found in the papers we cite when introducing each model.

Consider a model with compartments $X_1,X_2,\dots,X_n$. Then, the DFE has coordinates $X_i = 0$ for all $i \in \mathcal{I}$, where $\mathcal{I}$ is the set of the indexes of infectious compartments, and the EE, which we indicate with $(X_1^*,X_2^*,\dots,X_n^*)$, has all positive entries. A vast majority of Lyapunov functions in epidemic modelling fall into one of the categories listed below.
\begin{enumerate}
    \item \textbf{Linear combination of infectious compartments}. The Lyapunov function for the DFE when $\R0<1$ is of the form
    \begin{equation}
        L=\sum_{i \geq 2} c_i X_i,
    \end{equation}
    for some constants $c_i\geq 0$ to be determined \cite{bichara2014global,guihua2004global,guo2012global,hethcote2000mathematics,korobeinikov2004lyapunov,li2005global,li2002global,melesse2010global,ojo2017lyapunov,oke2019mathematical,syafruddin2013lyapunov,van1999global,yang2010global}. To prove convergence of the system to the DFE in this case it is often required the use of additional tools, such as LaSalle's invariance principle, which we briefly recall at the end of Section \ref{sec:SIS}. 
    \item \textbf{Goh-Lotka-Volterra}. The Lyapunov function for the EE when $\R0>1$ is of the form
    \begin{equation}
        L=\sum_i c_i (X_i - X_i^* \ln{X_i}),
    \end{equation}
    for some constants $c_i\geq 0$ to be determined \cite{ansumali2020modelling,bichara2012global,bichara2014global,harko2014exact,jardon2021geometric,kelkile2018stability,korobeinikov2004lyapunov,korobeinikov2002lyapunov,melesse2010global,ojo2017lyapunov,ottaviano2022globalmulti,ottaviano2022global,syafruddin2013lyapunov,tiantian2013global,vargas2009constructions}. These functions are adapted from a first integral of the notorious Lotka-Volterra prey-predator system, and were popularized by Bean-San Goh in a series of paper \cite{goh1976global,goh1977global,goh1979stability}.
    \item \textbf{Quadratic}. The Lyapunov function for the EE when $\R0>1$ is of the common form
    \begin{equation}
        L=\sum_i c_i (X_i - X_i^*)^2,
    \end{equation}
    for some constants $c_i\geq 0$ to be determined, or the composite form
    \begin{equation}
        L= \left(\sum_i X_i - X_i^*\right)^2.
    \end{equation}
    Some examples can be found in    \cite{liao2012global,mabotsa2022mathematical,taneco2020stability,vargas2009constructions,vargas2011global}.
    \item \textbf{Integral Lyapunov}. Lyapunov functions given as integrals over the dynamics of the model. The integration interval often start at some EE value $X_i^*$ and ends at the same $X_i$; this construction is very convenient if uniqueness of the EE is guaranteed, but the exact values of the EE are hard (or impossible) to determine analytically. Integral Lyapunov functions are particularly useful when the model includes multiple stages of infection, and consequently the infectious period changes from an exponential distribution to a gamma distribution \cite{cheng1982some,georgescu2013lyapunov,guo2012global,li2016class,sun2011global,tang2017new}. Integral Lyapunov functions, albeit in different forms, are widely used in models which incorporate explicit delay, such as systems of Delay Differential Equations (from here onwards, DDEs), and age-structured models. However, these fall beyond the scope of this paper, and we will briefly comment on them in Section \ref{Conc}.
    \item \textbf{Hybrid}. A linear combination of the above, which often includes the Goh-Lotka-Volterra in at least a few of the compartments of the system \cite{gonzalez2021qualitative,jardon2021geometric,li2021new,meskaf2020global,oke2019mathematical,ottaviano2022global,o2010lyapunov,tiantian2013global}.
\end{enumerate}

For some high-dimensional models, proving convergence to the EE might require additional tools, such as the geometric approach used in \cite{ottaviano2022global,van1999global}.

Lastly, we must notice that not all compartmental models only exhibit convergence to equilibrium. Some systems of autonomous ODEs may present stable or unstable limit cycles \cite{dafilis2012influence,ruan2003dynamical,wang2006epidemic}, homoclinic orbits \cite{ruan2003dynamical} or even chaos \cite{stiefs2009evidence}. In such cases, clearly, no global Lyapunov function may exist.

In the remainder of this survey, we will present various models and the corresponding Lyapunov functions, covering all the cases listed above.

\section{Epidemic models}

In this section, we present various compartmental epidemic models with the corresponding Lyapunov function(s).  We present the models from the smallest to the largest, in terms of number of compartments. We refer to \cite{anderson2013population,keeling2008infectious} for a basic introduction on compartmental epidemic models, and to \cite{shuai2013global} for a detailed exemplification of Lyapunov theory in this setting.

We provide a schematic representation of the flows in most of the systems we present. Flow diagrams can be useful to provide a visual, intuitive interpretation of the parameters involved in each system. Arrows between compartments indicate a change in the current state of individuals with respect to the ongoing epidemics, whereas arrows inward/outward the union of the compartments represent birth rate and death rate in the population. Often, these last two rates are considered to be equal, as this assumption allows the population to either remain constant or converge to a constant value, reducing the dimensionality of the system and (hopefully) its analytical complexity. However, some models include additional disease-induced mortality, to increase realism when modelling severe infectious diseases. We uniform the notation throughout the various models we present in this survey as much as possible, and provide a brief description of each parameter the first time it is encountered. We remark that each variable is assumed to be non-negative, since it represents a fraction of the population, but the biologically relevant region varies depending on the specific model we are describing.

Moreover, we illustrate the corresponding Lyapunov functions for 2D models, showcasing a selection of their power levels. The same procedure can be easily adapted to 3D models, but the corresponding visualizations can be hard to interpret in a static image.

\subsection{SIS}\label{sec:SIS}
The SIS model is characterized by the total absence of immunity after infection, i.e. the recovery from infection is followed by an instantaneous return to the susceptible class. The ODEs system which describes this situation is

\begin{minipage}{0.5\textwidth}
\begin{equation}\label{SIS}
    \begin{split}
        \dfrac{\text{d}S}{\text{d}t} &= \gamma I - \beta \dfrac{S I}{N},\\
        \dfrac{\text{d}I}{\text{d}t} &= \beta \dfrac{S I}{N} - \gamma I,
    \end{split}
\end{equation}
\vspace{.1cm}
\hfill
\end{minipage}
\begin{minipage}{0.6\textwidth}
\begin{tikzpicture}
\node[comp] (S) at (0,0) {$S$};
\node[comp] (I) at (2.5,0) {$I$};
\draw [-latex] (S) to node[above] {$\beta\frac{SI}{N}$} (I);
\draw [bend left, -latex] (I) to node[below] {$\gamma I$}  (S);
\end{tikzpicture}
\end{minipage}
where $\beta$ is the transmission rate and $\gamma$ is the recovery rate.

Notice that the population $N=S+I$ is constant, thus we can normalize it to $N=1$. Moreover, since $S+I=1$, we can reduce the system to one ODE which involves only infectious individuals
$$
\dfrac{\text{d}I}{\text{d}t} = (\beta (1-I) - \gamma) I.
$$
System \eqref{SIS} always admits the DFE, i.e. $E_0 = (1,0)$, and the EE, i.e. $E^* = \left(\dfrac{\gamma}{\beta},\dfrac{\beta - \gamma}{\beta}\right) $, which exists if and only if $\beta > \gamma$ (or equivalently if $\R0=\beta/\gamma>1$). Notice that, if $\R0<1$, then $I$ is always decreasing in the biologically relevant interval $[0,1]$.

A variation of model \eqref{SIS} can be obtained by adding demography to the system. This is the example of \cite{vargas2009constructions}, in which the authors consider a birth/immigration rate different from the natural death rate; moreover, they include an additional disease-induced death rate from infectious class. Thus, the population is not constant and the system of ODEs which describe the model is

\begin{figure}[H]
        \centering
        \begin{subfigure}[b]{0.49\textwidth}
            \centering
            \includegraphics[width=\textwidth]{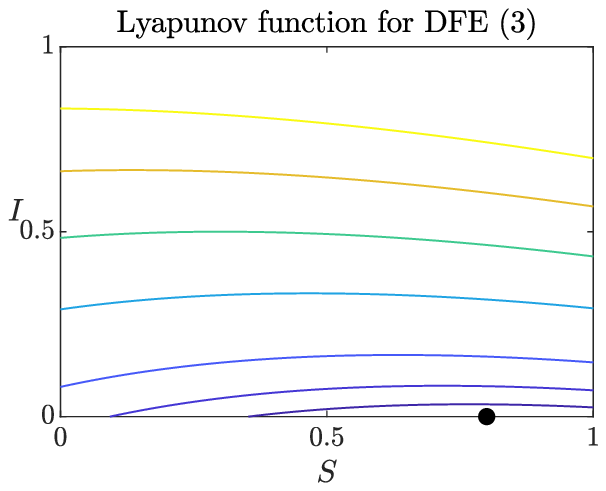}
            \caption{}%
            \label{fig:DFE3}
        \end{subfigure}
        \hfill
        \begin{subfigure}[b]{0.49\textwidth}  
            \centering 
            \includegraphics[width=\textwidth]{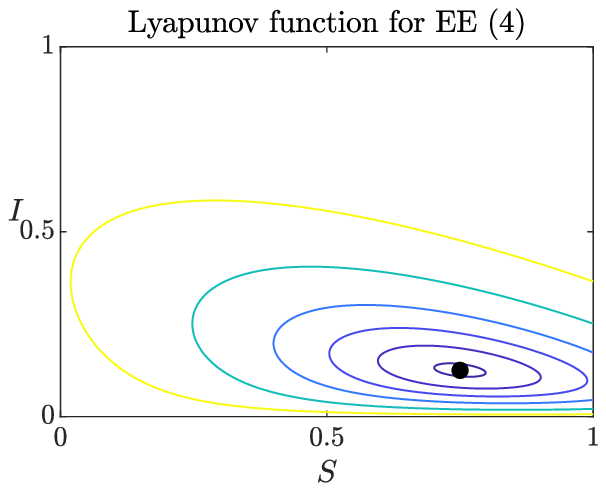}
            \caption{}%
            \label{fig:EE4}
        \end{subfigure}
        \vskip\baselineskip
        \begin{subfigure}[b]{0.49\textwidth}   
            \centering 
            \includegraphics[width=\textwidth]{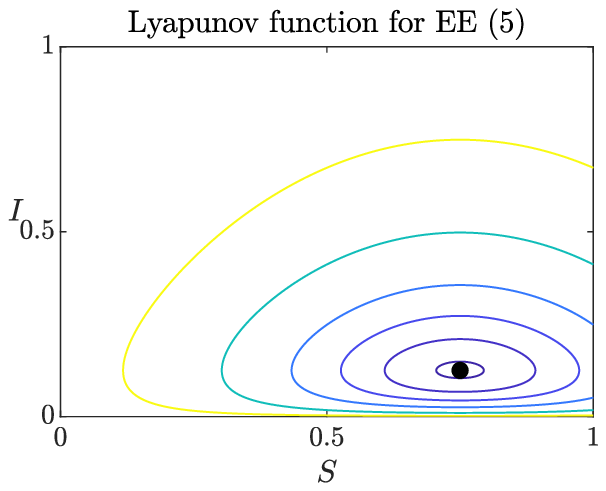}
            \caption{}%
            \label{fig:EE5}
        \end{subfigure}
        \hfill
        \begin{subfigure}[b]{0.49\textwidth}   
            \centering 
            \includegraphics[width=\textwidth]{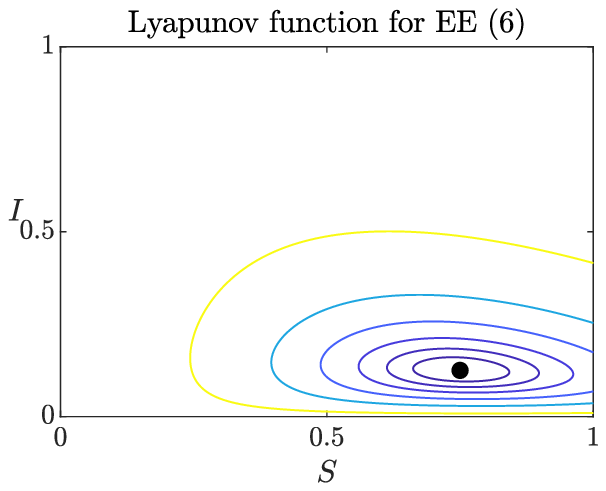}
            \caption{}%
            \label{fig:EE6}
        \end{subfigure}
        \caption{Power levels of Lyapunov functions \eqref{sis_DFE} (a), \eqref{sis_EE1} (b), \eqref{sis_EE2} (c), and \eqref{sis_EE3} (d). Values of the parameters are $\Lambda=0.8$, $\mu=1$, $\delta=1$, $\gamma=1$ in all the figures, $\beta=1$ in (a), so that $\R0=4/15<1$, and $\beta=4$ in (b), (c) and (d), so that $\R0=16/15>1$. We represent $V(S,I)=k$, with $k\in \{ 0.1, 0.25, 0.5, 1, 1.5, 2, 2.5 \}$ in (a), $k\in \{ 0.001, 0.01, 0.025, 0.05, 0.1, 0.2 \}$ in (b) and (c), and $k\in \{ 0.01, 0.025, 0.05, 0.1, 0.2, 0.5 \}$ in (d). Black dots represent the globally stable equilibrium the system converges to, and correspond to $V(S,I)=0$.} 
        \label{fig:2DLyap}
    \end{figure}

\begin{minipage}{0.5\textwidth}
\begin{equation}\label{SIS_demo}
    \begin{split}
        \dfrac{\text{d}S}{\text{d}t} &= \Lambda + \gamma I - \beta \dfrac{S I}{N} - \mu S,\\
        \dfrac{\text{d}I}{\text{d}t} &= \beta \dfrac{S I}{N} - (\delta +\gamma+\mu) I,
    \end{split}
\end{equation}
\hfill
\end{minipage}
\begin{minipage}{0.6\textwidth}
\begin{tikzpicture}
\node[comp] (S) at (0,0) {$S$};
\node[comp] (I) at (2.5,0) {$I$};
\node (void) at (-1.5,0) {};
\node (void2) at (0,-1.5) {};
\node (void3) at (2.5,-1.5) {};
\draw [-latex] (S) to node[above] {$\beta\frac{SI}{N}$} (I);
\draw [bend left, -latex] (I) to node[below] {$\gamma I$}  (S);
\draw [-latex] (void) to node[above] {$\Lambda$} (S);
\draw [-latex] (S) to node[left] {$\mu S\,$} (void2);
\draw [-latex] (I) to node[right] {$(\delta + \mu)I$} (void3);
\end{tikzpicture}
\vspace{.1cm}
\end{minipage}
\medskip

\noindent
where $\Lambda$ represents the birth/immigration rate, $\mu$ the natural death rate and $\delta$ the disease-induced mortality rate. System \eqref{SIS_demo} always admits the DFE, namely $E_0=(S_0,0):=\left(\dfrac{\Lambda}{\mu},0\right)$, and the EE, namely $E^*=(S^*,I^*)$, where $I^*>0$ if and only if $\R0 = \dfrac{\Lambda \beta}{\mu (\mu + \delta + \gamma)} > 1$. 
In \cite{vargas2009constructions}, a Lyapunov function for the DFE is defined as
\begin{equation}\label{sis_DFE}
V(S,I):= \dfrac{1}{2}\left(S-S_0 + I\right)^2 + \dfrac{2 \mu + \delta}{\beta} I,
\end{equation}
whereas the Lyapunov function for the EE is built using a combination of the quadratic and logarithmic functions
\begin{equation}\label{sis_EE1}
    V(S,I):= \dfrac{1}{2}\left(S-S^* + I-I^*\right)^2 + \dfrac{ 2 \mu + \delta }{\beta}\left(I - I^* -I^*\ln\left(\dfrac{I}{I^*}\right)\right).
\end{equation}
The authors also construct two more examples of Lyapunov functions for the EE, namely
\begin{equation}\label{sis_EE2}
V(S,I):= \dfrac{1}{2}(S-S^*)^2 + \dfrac{\mu + \delta}{ \beta} \left( I- I^*- I^* \ln\left(\dfrac{I}{I^*}\right)\right),
\end{equation}
and
\begin{equation}\label{sis_EE3}
\begin{split}
    V(S,I):= &\dfrac{1}{2}\left( S-S^* +I-I^*\right)^2 + \dfrac{S^*(\delta + 2 \mu)}{2 \gamma} \left( S- S^*- S^* \ln\left(\dfrac{S}{S^*}\right)\right) \\
    &+ \dfrac{S^* (\delta + 2\mu)}{\gamma} \left( I- I^*- I^* \ln\left(\dfrac{I}{I^*}\right)\right).
\end{split}
\end{equation}
Power levels of the functions \eqref{sis_DFE}, \eqref{sis_EE1}, \eqref{sis_EE2} and \eqref{sis_EE3} are visualized if Figure \ref{fig:2DLyap}. By definition of a Lyapunov functions, orbits of the corresponding system \eqref{SIS_demo} evolve on decreasing power levels, and they tend to the corresponding equilibrium as $t \to +\infty$.
In \cite{vargas2011global} the author found a simpler Lyapunov function for the DFE when $\R0<1$, i.e.
\begin{equation}\label{quadratic}
    V(I) = \dfrac{1}{2}I^2.
\end{equation}
However, this last Lyapunov function \eqref{quadratic} only ensures that $I \to 0$ as $t \to +\infty$. To complete the proof of the converge of the system to the DFE, one needs in addiction to invoke LaSalle's theorem \cite{lasalle1976stability} (see also \cite[Thm. 3.4]{khalil2015nonlinear}), as is indeed done in \cite{vargas2011global}.

Considering the importance of this theorem, especially when combined with the use of Lyapunov functions, we include its statement here.
\begin{theorem} (LaSalle's invariance principle) Let $X'=f(X)$ be a system of $n$ ODEs defined on a positively invariant set $\Omega\subset\mathbb{R}^n$. Assume the existence of a function $V\in C^1(\Omega,\mathbb{R})$ such that $V'(X)\leq 0$ for all $X \in \Omega$. Let $M_V$ be the set of stationary points for $V$, i.e. $V'(X)=0$ for all $X \in M_V$, and let $N$ be the largest invariant set of $M_V$. Then, every solution which starts in $\Omega$ approaches $N$ as $t\to +\infty$.
\end{theorem}
In particular, this theorem implies that, if we can prove the approach of the disease to the manifold describing absence of infection and the uniqueness of the DFE, then the DFE is GAS.

\subsection{SIR/SIRS}

The SIR model is characterized by the total immunity after the infections, i.e. recovered individuals can not become susceptible again. A classical example for this scenario is measles. The ODEs system which describes this situation is

\begin{minipage}{0.5\textwidth}
\begin{equation}\label{SIR}
    \begin{split}
        \dfrac{\text{d}S}{\text{d}t} &= - \beta \dfrac{S I}{N},\\
        \dfrac{\text{d}I}{\text{d}t} &= \beta \dfrac{S I}{N} - \gamma I,\\
        \dfrac{\text{d}R}{\text{d}t} &= \gamma I,
    \end{split}
\end{equation}
\vspace{.1cm}
\hfill
\end{minipage}
\begin{minipage}{0.6\textwidth}
\begin{tikzpicture}
\node[comp] (S) at (0,0) {$S$};
\node[comp] (I) at (2,0) {$I$};
\node[comp] (R) at (2,-2) {$R$};
\draw [-latex] (S) to node[above] {$\beta\frac{SI}{N}$} (I);
\draw [-latex] (I) to node[midway, left] {$\gamma I$}  (R);
\end{tikzpicture}
\end{minipage}
\medskip

\noindent
where $\beta$ is the transmission rate and $\gamma$ is the recovery rate. 

If we assume that recovered individuals eventually lose their immunity, we obtain the SIRS model. Denoting by $\alpha$ the immunity loss rate, we obtain the following ODEs system

\begin{minipage}{0.5\textwidth}
\begin{equation}\label{SIRS}
    \begin{split}
        \dfrac{\text{d}S}{\text{d}t} &= - \beta \dfrac{S I}{N} + \alpha R,\\
        \dfrac{\text{d}I}{\text{d}t} &= \beta \dfrac{S I}{N} - \gamma I,\\
        \dfrac{\text{d}R}{\text{d}t} &= \gamma I - \alpha R.
    \end{split}
\end{equation}
\vspace{.1cm}
\hfill
\end{minipage}
\begin{minipage}{0.6\textwidth}
\begin{tikzpicture}
\node[comp] (S) at (0,0) {$S$};
\node[comp] (I) at (2,0) {$I$};
\node[comp] (R) at (2,-2) {$R$};
\draw [-latex] (S) to node[above] {$\beta\frac{SI}{N}$} (I);
\draw [-latex] (I) to node[midway, right] {$\gamma I$}  (R);
\draw [-latex] (R) to node[midway, left] {$\alpha R $} (S);
\end{tikzpicture}
\end{minipage}
It is clear that, if $\alpha = 0$, system \eqref{SIRS} coincides with system \eqref{SIR}.

These models admit only the DFE; in order to have an EE, we need to add the demography to model \eqref{SIR} or \eqref{SIRS}. 

In \cite{vargas2009constructions}, the authors consider the following ODEs system

\begin{minipage}{0.5\textwidth}
\begin{equation}\label{SIRS_demo}
    \begin{split}
        \dfrac{\text{d}S}{\text{d}t} &= \Lambda - \beta \dfrac{S I}{N} -\mu S + \alpha R,\\
        \dfrac{\text{d}I}{\text{d}t} &= \beta \dfrac{S I}{N} -( \gamma+\delta+\mu) I,\\
        \dfrac{\text{d}R}{\text{d}t} &= \gamma I - (\alpha+\mu)R.
    \end{split}
\end{equation}
\hfill
\end{minipage}
\begin{minipage}{0.45\textwidth}
\begin{tikzpicture}
\node[comp] (S) at (0,0) {$S$};
\node[comp] (I) at (2.5,0) {$I$};
\node[comp] (R) at (2.5,-2) {$R$};
\node (void) at (-1.5,0) {};
\node (void2) at (0,-2) {};
\node (void3) at (4,0) {};
\draw [-latex] (S) to node[above] {$\beta\frac{SI}{N}$} (I);
\draw [-latex] (I) to node[midway, right] {$\gamma I$}  (R);
\draw [-latex] (R) to node[midway, left] {$\alpha R\,$} (S);
\draw [-latex] (void) to node[above] {$\Lambda$} (S);
\draw [-latex] (S) to node[left] {$\mu S\,$} (void2);
\draw [-latex] (R) to node[below] {$\mu R$} (void2);
\draw [-latex] (I) to node[above] {$(\delta+\mu)I$} (void3);
\end{tikzpicture}
\vspace{.1cm}
\end{minipage}
System \eqref{SIRS_demo} admits the DFE, $E_0=(S_0,0,0)$, and the EE, $E^*=(S^*,I^*,R^*)$, which exists if and only if $\R0 =\dfrac{\beta \Lambda}{\mu(\mu + \gamma + \delta) } >1$. In \cite{vargas2009constructions}, the Lyapunov function for the DFE is defined as follows
$$
V(S,I,R) := \dfrac{1}{2}\left( S-S_0 + I + R\right)^2 + \dfrac{2 \mu + \delta }{\beta }I + \dfrac{2 \mu + \delta}{2 \gamma} R^2,
$$
whereas the Lyapunov function for the EE is the combination of the composite quadratic, common quadratic and logarithmic functions as follows
\begin{equation*}
    \begin{split}
        V(S,I,R) := & \dfrac{1}{2}\left(S-S^* + I-I^* + R-R^* \right)^2  \\
        & + \dfrac{2 \mu + \delta}{\beta} \left( I-I^* - I^* \ln\left( \dfrac{I}{I^*}\right)\right) + \dfrac{2\mu + \delta}{2 \gamma} (R-R^*)^2.
    \end{split}
\end{equation*}
The authors also present other Lyapunov functions for SIR/SIRS models; in particular, they also cite \cite{beretta2001global,mena1992dynamic} in which some variations of system \eqref{SIRS_demo} are showed. Other Lyapunov functions for SIR/SIRS epidemic models are in \cite{shuai2013global}, in which the authors use a graph-theoretic approach. 

In \cite{vargas2011global}, the author proved that the quadratic Lyapunov function \eqref{quadratic} of the SIS model applies to the SIR and the SIRS, as well.

\subsection{SEIR/SEIS/SEIRS}

In \cite{korobeinikov2004lyapunov}, the authors study both SEIR and SEIS models. Many real world examples present a phase of exposition to the disease, between susceptibility and infectiousness. The models presented thus far, albeit simpler to study, are unable to replicate this mechanism.

The authors first analyze a SEIR model with demography and constant population, in which the disease is transmitted both horizontally and vertically. Individuals infected vertically pass first in the exposed compartment. The ODEs system which describe the model is\\
\hspace{-.5cm}
\begin{minipage}{0.57\textwidth}
\begin{equation}\label{eq:SEIR}
    \begin{split}
        \dfrac{\text{d}S}{\text{d}t} =& \mu - \beta S I - p \mu I - q \mu E - \mu S,\\
        \dfrac{\text{d}E}{\text{d}t} =& \beta S I + p \mu I - \theta E - \mu E + q \mu E,\\
        \dfrac{\text{d}I}{\text{d}t} =& \theta E - (\delta + \mu)I,
    \end{split}
\end{equation}
\end{minipage}
\hspace{-1cm}
\begin{minipage}{0.4\textwidth}
\begin{tikzpicture}
\node[comp] (S) at (0,0) {$S$};
\node[comp] (E) at (2.5,0) {$E$};
\node[comp] (I) at (1,-1.5) {$I$};
\node (void) at (-2,0) {};
\node (void2) at (0,2) {};
\node (void3) at (2.5,2) {};
\node (void4) at (-1,-1.5) {};
\node (void5) at (2.5,-1.5) {};
\draw [-latex] (S) to node[above] {$\beta S I$} (E);
\draw [-latex] (E) to node[above] {$\theta E\ \ $}  (I);
\draw [-latex] (void2) to node[left] {$\mu(1-pI-qE)$} (S);
\draw [-latex] (S) to node[above] {$\mu S$} (void);
\draw [-latex] (void3) to node[left] {$\mu(pI+qE)$} (E);
\draw [-latex] (E) to node[right] {$\mu E$} (void5);
\draw [-latex] (I) to node[above] {$(\delta + \mu)I$} (void4);
\end{tikzpicture}
\vspace{.1cm}
\end{minipage}
\medskip

\noindent
and $R= 1- S-E-I$. The vertical transmission of the disease is represented by the probabilities $p$ and $q$ of being born directly in the Exposed compartment, rather than in the Susceptible one, and is represented by the inward arrow in compartment E.

The authors first provide an equivalent system, performing the substitution $(S,E,I) \longrightarrow (P,E,I)$, where $P:=S + p \dfrac{\mu}{\beta}$. They then proceed to prove the GAS of the EE, using the following Lyapunov function
\begin{equation*}
    \begin{split}
V(P,E,I) := &(P-P^* \ln P) + \dfrac{\theta + \mu}{\theta + \mu - q \mu } (E-E^* \ln E) \\&+ \dfrac{\theta +\mu }{\theta + \mu - q \mu } ( I- I^* \ln I).
    \end{split}
\end{equation*}

Later, the authors analyze a situation in which the recovery does not provide immunity, namely the SEIS model. They also assume that a fraction $r$ of offspring of the infective hosts is born directly into the infective compartment. In this case, the ODEs system changes accordingly
describe the model is
\begin{equation}\label{eq:SEIS}
    \begin{split}
        \dfrac{\text{d}S}{\text{d}t} =& \mu - \beta S I + (\delta - p \mu - r \mu) I - q \mu E - \mu S,\\
        \dfrac{\text{d}E}{\text{d}t} =& \beta S I + p \mu I - (\theta + \mu - q \mu )E,\\
        \dfrac{\text{d}I}{\text{d}t} =& \theta E - (\delta + \mu - \mu r)I,
    \end{split}
\end{equation}
and $S+E+I = 1$. Notice that, due to the population remaining constant in system \eqref{eq:SEIS}, one could in principle reduce its dimensionality and consider it as a planar system.

The authors prove the GAS of the EE using the following Lyapunov function
\begin{equation*}
    \begin{split}
V(S,E,I) :=& (S-S^* \ln S) + \mu \dfrac{1-S^*}{\beta I^*S^* } (E-E^* \ln E)\\& + \mu \dfrac{1-S^*}{\theta E^* }\left( 1 + p \rho_0 \dfrac{\mu}{\beta}\right)   ( I- I^* \ln I).
    \end{split}
\end{equation*}
A natural extension to these models is the SEIRS \cite{hethcote1991some,van1999global}, in which one can combine the existence of an immune compartment and the loss of immunity. It is described by the following system of ODEs\\
\begin{minipage}{0.55\textwidth}
\begin{equation}\label{eq:SEIRS}
    \begin{split}
        \dfrac{\text{d}S}{\text{d}t} =& - \beta g(I) S + \mu -\mu S + \alpha R ,\\
        \dfrac{\text{d}E}{\text{d}t} =& \beta g(I) S  - (\theta + \mu )E,\\
        \dfrac{\text{d}I}{\text{d}t} =& \theta E - (\gamma + \mu)I,\\
        \dfrac{\text{d}R}{\text{d}t} =& \gamma I - (\alpha + \mu)R,
    \end{split}
\end{equation}
\hfill
\end{minipage}
\begin{minipage}{0.4\textwidth}
\begin{tikzpicture}
\node[comp] (S) at (0,0) {$S$};
\node[comp] (E) at (2.5,0) {$E$};
\node[comp] (I) at (2.5,-2.5) {$I$};
\node[comp] (R) at (0,-2.5) {$R$};
\node (void) at (0,1.5) {};
\node (void2) at (-1,-1.25) {};
\node (void3) at (3.5,-1.25) {};
\draw [-latex] (S) to node[above] {$\beta g(I) S$} (E);
\draw [-latex] (I) to node[above] {$\gamma I$}  (R);
\draw [-latex] (R) to node[right] {$\alpha R\,$} (S);
\draw [-latex] (E) to node[left] {$\theta E$} (I);
\draw [-latex] (E) to node[above] {$\ \ \ \mu E$} (void3);
\draw [-latex] (void) to node[left] {$\mu$} (S);
\draw [-latex] (S) to node[above] {$\mu S\ \,$} (void2);
\draw [-latex] (R) to node[below] {$\mu R\ \ \,$} (void2);
\draw [-latex] (I) to node[below] {$\ \mu I$} (void3);
\end{tikzpicture}
\vspace{.1cm}
\end{minipage}
\medskip

\noindent
where $g\in C^3(0,1]$, $g(0)=0$ (meaning, in absence of infectious individuals, the disease does not spread) and $g(I)>0$ for $I>0$. The classical choice is $g(I)=I$, as in systems \eqref{eq:SEIR} and \eqref{eq:SEIS}. 
Assuming moreover
$$
\lim_{I \to 0^+} \dfrac{g(I)}{I}=c\in [0,+\infty),
$$
the authors of \cite{hethcote1991some} derive $\R0=\dfrac{c\beta\theta}{(\theta+\mu)(\gamma+\mu)}$. They then prove GAS of the DFE of system \eqref{eq:SEIRS} through the use of the following linear Lyapunov function
$$
V(E,I)=E+\dfrac{\theta+\mu}{\theta}I,
$$
whereas the GAS of the EE is proved with a more complex geometrical method in \cite{van1999global}.

\subsection{SAIR/SAIRS}

One of the main challenges of the Covid-19 pandemic was the presence of asymptomatic individuals spreading the disease. Such individuals must clearly be somehow distinguished from symptomatic infectious individuals, as they are likely to behave like a susceptible individual. Even though their viral load, and hence infectiousness, might be smaller, they are more likely to get in close contact with susceptible individuals.

In \cite{ottaviano2022global}, the authors consider a SAIRS model. The main difference between this kind of models and the SEIR is that both asymptomatic and symptomatic hosts may infect susceptible individuals. The immunity is not permanent, i.e. recovered individuals will become susceptible again after a certain period of time. Moreover, vaccination are included. The ODEs system which describe this model is

\begin{minipage}{0.55\textwidth}
\begin{equation} \label{sairs_s}
\begin{split}
     \frac{\text{d} S}{\text{d}t} &= \mu  - \bigg(\beta_A A + \beta_I I\bigg)S -(\mu + \nu) S +\gamma R,\\
     \frac{\text{d} A}{\text{d}t} &=  \bigg(\beta_A A + \beta_I I\bigg)S -(\alpha + \delta_A +\mu) A, \\ 
     \frac{\text{d} I}{\text{d}t} &= \alpha A - (\delta_I + \mu)I, \\ 
     \frac{\text{d} R}{\text{d}t} &=  \delta_A A +\delta_I I + \nu S - (\gamma + \mu)R,
     \end{split}
\end{equation}
\hfill
\end{minipage}
\begin{minipage}{0.5\textwidth}
\begin{tikzpicture}[scale=1]
\node (void) at (0,1.5) {};
\node (void2) at (2,1.5) {};
\node (void3) at (2,1.5) {};
\node (void4) at (2,-4.5) {};
\node (void5) at (2,-4.5) {};
\node[comp] (S) at (0,0) {$S$};
\node[comp] (A) at (4,0) {$A$};
\node[comp] (R) at (0,-3) {$R$};
\node[comp] (I) at (4, -3) {$I$};
\draw[-latex] (void) to node[left] {$\,\mu$} (S);
\draw[-latex] (S) to node[left] {$\mu S$} (void2);
\draw[-latex] (S) to node[above] {$(\beta_A A + \beta_I I)S $} (A);
\draw[-latex] (I) to node[above] {$\delta_I I$}  (R);
\draw[-latex] (R) to node[right] {$\gamma R\ $} (S);
\draw[-latex] (A) to node[left] {$\delta_A A\,$} (R); 
\draw[-latex] (A) to node[right] {$\ \mu A$} (void3);
\draw[-latex] (A) to node[right] {$\alpha A\,$} (I);
\draw[latex-] (R) edge[bend left] node[left] {$\nu S$} (S.south west);
\draw[-latex] (I) to node[right] {$\mu I$} (void4);
\draw[-latex] (R) to node[left] {$\mu R\ $} (void5);
\end{tikzpicture}
\end{minipage}
\medskip

\noindent
The global stability analysis of the EE has been performed for two variations of the original model, described in the following.

The first model analyzed is the SAIR model, i.e. the case in which recovery from the disease grants permanent immunity. In this case, the corresponding Lyapunov function is the combination of the Lokta-Volterra Lyapunov functions for $S$, $A$ and $I$
\begin{equation*}
\begin{split}
V(S,A,I) :=& c_1 S^* \left( \dfrac{S}{S^*}- 1 - \ln \left( \dfrac{S}{S^*} \right) \right) + c_2  A^* \left( \dfrac{A}{A^*}- 1 - \ln \left( \dfrac{A}{A^*} \right) \right) \\
&+  I^* \left( \dfrac{I}{I^*}- 1 - \ln \left( \dfrac{I}{I^*} \right) \right)  ,
\end{split}
\end{equation*}
where $c_1, c_2 >0$.

The second model is the SAIRS model, with homogeneous disease transmission and recovery among $A$ and $I$, i.e. $\beta_A=\beta_I$ and $\delta_A = \delta_I$. In this case, it is possible to sum equations for $A$ and $I$, defining $M:= A + I$, reducing the dimensionality of the system. Thus, the Lyapunov function can be written as the combination of the square function and the Lokta-Volterra as follows
    $$
    V(S,M) := \dfrac{1}{2}(S-S^*)^2 + w \left( M-M^* -M^* \ln\left( \dfrac{M}{M^*}\right)\right),
    $$
    where $w >0$.

The global stability in the most general case is proved similarly to \cite{van1999global}.

\subsection{More exotic compartmental models}

The aforementioned models are some of the most commonly used in literature. In order to capture additional disease-specific nuances, these model can be modified or extended by adding new compartments.

Some diseases, for example, present different stages of infection. In this case, an infected individual can progress between two or more stages before recovering. In \cite{guo2012global}, the authors perform the global stability analysis via an integral Lyapunov function of a general class of multistage models. In their model, infectious individual can move both forward and backward on the chain of stages, in order to incorporate both a natural disease progression and the amelioration due to the effects of treatments.

The system of ODEs which describes the model is
\begin{equation}
    \begin{split}
         \frac{\text{d} S}{\text{d}t} &= \theta(S) - f(N) \sum_{j=1}^n g_j(S,I_j), \\
          \frac{\text{d} I_1}{\text{d}t} &= f(N) \sum_{j=1}^n g_j(S,I_j) + \sum_{j=1}^n \phi_{1,j}(I_j) - \sum_{j=1}^{n+1}\phi_{j,1}(I_1)-\zeta_1(I_1),\\
          \frac{\text{d} I_i}{\text{d}t} &= \sum_{j=1}^n \phi_{i,j}(I_j) - \sum_{j=1}^{n+1}\phi_{j,i}(I_i)-\zeta_i(I_i),\qquad i=2,3,\dots,n,
    \end{split}
\end{equation}
where $\theta(S)$ is the growth function, 
 $f(N)\sum_{j=1}^n g_j(S,I_j)$ is the incidence term, $\zeta_i(I_i)$, $1\leq i \leq n$, denote the removal rates of the $I_i$ compartment. Moreover, for any $i,j=1,\dots,n$, the functions $\phi_{i,j}(I_j)$ represent the rate of the disease progression if $i>j$ and the amelioration if $i < j$.

The corresponding Lyapunov function for the DFE is linear in the disease compartments, i.e.
\begin{equation}
V(I_1,\dots,I_n) = \sum_{i=1}^n c_i I_i,     
\end{equation}
where $c_1 = \mathcal{R}_0$ and $c_i \geq 0$ for all $i=2,\dots,n$. For the global stability of the EE the authors made some assumptions on the aforementioned functions. In particular, they consider the following integral Lyapunov function
\begin{equation}
V(S,I_1,\dots,I_n) = \tau \int_{S^*}^{S} \dfrac{\Phi(\xi)-\Phi(S^*)}{\Phi(\xi)} \ \text{d}\xi + \sum_{i=1}^n \tau_i \int_{I_i^*}^{I_i} \dfrac{\psi_i(\xi) - \psi_i(I_i^*)}{\psi_i(\xi)} \ \text{d}\xi,
\end{equation}
where $\tau, \tau_i > 0$, for all $i=1,\dots, n$. For a more in-depth explanation on the functions $\Phi(\cdot)$ and $\psi_i(\cdot)$ we refer to \cite[Sect.~5]{guo2012global}.

Diseases which present multiple virus strains, due to the existence of different serotypes of the virus or due to a mutation of the original disease, may need to be modelled differently. Dengue, tuberculosis and various sexually transmitted diseases are caused by more than one strain of a pathogen. Influenza type A viruses mutate constantly: an infection with one of its strains gives permanent immunity against that specific strain. However, the so called ``antigenic drift'' produces new virus strains, thus the hosts only acquire partial immunity, or no immunity at all. Modelling these types of diseases requires the inclusion of cross-protective effects, in which the immunity acquired towards one strain offers partial protection towards another strain based on their antigenic similarity. In \cite{bichara2014global}, the authors consider an $n$ strain model, both without immunity and with immunity for all the strains. Moreover, they analyze an MSIR model, in which the $M$ compartment represents the proportion of newborns who possess temporary passive immunity due to protection from maternal antibodies.
For all the three model, the authors use a linear Lyapunov function to prove the global stability of the DFE and a logarithmic Lyapunov function to prove the global stability of the EE.

Other compartmental models include e.g. control strategies. For new ongoing epidemics, the most immediate strategy is including quarantine and isolation of infectious individuals. 
For well-known epidemics for which a vaccination is available, it is useful to incorporate a vaccinated individuals compartment $V$ to keep track of the two possible immunities, disease and vaccine induced, respectively. Usually, vaccination does not confer permanent immunity, and after a certain disease-dependent period individuals become susceptible again. An example is \cite{oke2019mathematical}, in which the authors analyze a SIRV epidemic model with non-linear incidence rate. The global stability of the DFE is proved using as linear Lyapunov function the infectious compatment $I$ and the global stability of the EE, instead, using a combination of a quadratic function in $S$ and a logarithmic function in the compartments $I$ and $V$.

\section{Conclusion}\label{Conc}

In this survey, we presented the most widely used Lyapunov functions in the field of epidemic compartmental models. We focused on systems expressed as autonomous systems of ODEs. These models allow for various interesting generalizations, of which we provide a non-comprehensive list below.

One extension of the classic compartmental epidemic models is the so-called multi-group approach, see e.g. \cite{kuniya2013global,sun2011global}. These models describe $n$ communities, interacting with each other, and whose internal evolution follows a standard compartmental model. 
A first example of such a model is presented in \cite{fall2007epidemiological}, in which the authors consider a $n$ groups SIS model. In order to prove the GAS of the EE, they use a results on Metzler matrices. In \cite{shuai2013global}, the authors consider a heterogeneous SIS disease model, for which they provide Lyapunov functions both for the DFE and for the EE. For the latter, they use a complex graph-teoretic method, for the details of which we refer to the original paper. Global stability of EE via Lyapunov function for multi-group generalization can be found also for the SIR \cite{guo2006global}, SIRS \cite{muroya2013global}, SEIR \cite{guo2008graph} and SAIR/SAIRS model \cite{ottaviano2022globalmulti}. Notice that, due to the complexity of the models, some of them require additional technical assumptions to prove the global stability of the endemic equilibrium. 

Other classes of models include interactions between human and vector population, i.e. animals which transmit the disease to humans, or with the pathogens, such as viruses or bacteria. In both cases, authors often include a compartmental structure for the non-human population. Some examples of vector-host models are shown in \cite{syafruddin2013lyapunov,tewa2009lyapunov,yang2010global}. Another example can be found in \cite{liao2012global}, in which a SIR-B compartmental model is considered. Here the ``B'' denotes the concentration of the pathogen in the environment. 

All the models discussed thus far are described by only autonomous systems of ODEs. However, in order to increase realism, it is possible to use non-autonomous systems to describe the spread of an infectious disease. This is the case of systems in which some parameters change in time \cite{martcheva2009non,sottile2020time}, to describe seasonal changes, or in which the state variables depend on the previous state, i.e. the model includes a time delay \cite{beretta1995global,wang2002global}. In these cases, it is still possible to find Lyapunov functions to prove the global stability of the equilibria using other techniques, described for example in \cite{lasalle1976stability}.

Another popular option is to explicitly include delay in the system, such as in \cite{beretta1995global,huang2013note,Huang2010viral,huang2010global,mccluskey2010complete, tiantian2013global,xu2009global}. In the latter the authors perform the global stability analysis of a SEIQR model, in which $Q$ denotes the quarantined individuals. They explicitly include a latent period for the infection, transforming two of the ODEs in DDEs. The corresponding Lyapunov function includes the integration over an interval whose size is precisely the latent period.

Lastly, a widely adopted strategy is to explicitly include the ``time since infection'' \cite{chekroun2019global,huang2012lyapunov,mccluskey2012global,yang2014global,yang2015global} in age-structured models. This allows to explicitly take into account time heterogeneity in the spread of an infectious disease in a population. These cases are outside of the scope of this project, and we leave them as inspiration for future works.

As a final remark, some recent results in a more theoretical approach to the topic are worth mentioning. They focus on existence and characteristics of such functions rather than on applications to epidemiological models. We refer the interested reader to \cite{abbondandolo2018chain,bernard2018lyapounov,fathi2015aubry,fathi2012smooth,hafstein2021smooth,paternain2002boundary} and the references therein.

\paragraph{Acknowledgments.}
The authors are grateful to the organizers of the conference \emph{100 Years Unione Matematica Italiana - 800 Years Università di Padova}, which made their scientific cooperation possible. Moreover, they acknowledge Politecnico di Milano, Polish Academy of Sciences, Inria and University of Trento for supporting their research.

\bibliographystyle{plain} 
\bibliography{biblio} 
\end{document}